\documentclass [a4 paper,12pt]{article}
\usepackage{mathrsfs}

\usepackage[usenames]{color}
\usepackage{bbm,amsmath,amssymb,amsfonts,amsthm,graphics,graphicx}

\newtheorem{lem}{Lemma}[section]
\newtheorem{thm}[lem]{Theorem}

\newtheorem{claim}{Claim}
\newtheorem{pro}{Proposition}
\newtheorem{con}{Conjecture}

\DeclareMathOperator{\Pro}{\mathbbm{P}}
\DeclareMathOperator{\V}{\mathbbm{Var}}
\DeclareMathOperator{\E}{\mathbbm{E}}
\DeclareMathOperator{\En}{\mathscr{E}}

\DeclareMathOperator{\si}{\sigma} \DeclareMathOperator{\la}{\lambda}
\DeclareMathOperator{\ep}{\epsilon}
\DeclareMathOperator{\rank}{rank}

\title{The limiting spectral distribution of the\\ generalized
Wigner matrix\footnote{Supported by NSFC No.10831001, PCSIRT and the
``973" program.}}

\author{ \small Wenxue Du,~~Xueliang Li,~~Yiyang Li \\
\small Center for Combinatorics and LPMC-TJKLC\\
\small Nankai University, Tianjin 300071, P.R. China\\
\small Email: lxl@nankai.edu.cn}

\date{}

\begin{document}
\maketitle

\begin{abstract}
The properties of eigenvalues of large dimensional random matrices
have received considerable attention. One important achievement is
the existence and identification of the limiting spectral
distribution of the empirical spectral distribution of eigenvalues
of Wigner matrix. In the present paper, we explore the limiting
spectral distribution for more general random matrices, and,
furthermore, give an application to the energy of general random
graphs, which generalizes the result of Nikiforov.\\[3mm]
{\bf Keywords}: eigenvalues, random matrix, Wigner matrix, empirical
spectral distribution, limiting spectral distribution, moment
approach, Stieltjes transform, graph energy.\\[3mm]
{\bf AMS Subject Classification 2000: 15A52, 15A18, 05C80, 05C90,
92E10}

\end{abstract}

\section{Introduction}

 In quantum mechanics, the energy levels of quanta can be
 characterized by the eigenvalues of a matrix. The empirical
 spectral distribution (ESD) of a matrix, however, is rather
 complicated when the order of the matrix is high. Wigner
 \cite{W55,W58} considered the limiting spectral distribution
 (LSD) for large dimensional random matrices, and obtained the
 famous semi-circle law. We recall a generalization here due
 to Bai \cite{B}. To be precise, Wigner investigated the
 LSD for a random matrix, so-called {\em Wigner matrix},
 $$\mathbf{X}_n:=(x_{ij}), ~~1\le i,j\le n,$$
 which satisfies the
following properties:\begin{itemize}

  \item $x_{ij}$'s are independent random variables with $x_{ij}=x_{ji}$;

  \item the $x_{ii}$'s have the same distribution $F_1$,
  while the $x_{ij}$'s are to
 possess
 the same distribution $F_2$;

 \item $\V(x_{ij})=\si_2^2<\infty$ for all $1\le i< j\le
 n$.
 \end{itemize}
 Set $$\mathbf{Y}_n=\frac{1}{2\sqrt{n}}\mathbf{X}_n.$$ We denote the eigenvalues of
 $\mathbf{Y}_n$ by $\la_{1,n},\la_{2,n},\ldots,\la_{n,n}$, and their
 ESD by $\Phi_n(x)=N_n(x)/n$ where
 $$N_n(x):=\#\{\la_{k,n}\mid\la_{k,n}\le x, ~k=1,2,\ldots, n\}.$$
 It is readily seen that for any given real number $x$, the ESD
 $\Phi_n(x)$ is a random variable on the space $\mathcal{X}_n$ consisting of Wigner
 matrices, while for any given
 matrix $\mathbf{X}_n$ in $\mathcal{X}_n$, $\Phi_n(x)$ can be regarded as a
 distribution function of the eigenvalues of $\mathbf{X}_n$.
\begin{thm}[\bf Wigner \cite{W55,W58}]\label{Thm-1} Let $F_1$ and $F_2$ be distribution
functions  mentioned above. Then
 $$\lim_{n\rightarrow\infty}
\Phi_n(x)=\Phi(x)\mbox{ a.s.,}$$
 where $\Phi(x)$ is the limiting spectral distribution
 with density
 $$\label{SC Law-0}
 \phi(x)=\left\{\begin{array}{ll}\displaystyle
  \frac{2}{\pi\si_2^2}\sqrt{\si_2^2- x^2 } & \mbox{if~~ } |x|\le\si_2,\\
  0 & \mbox{if~~ } |x|>\si_2. \end{array}\right.  $$
\end{thm}

 In the present paper, we explore the LSD for more general random
 matrices. Let $m:=m(n)\ge 2$ be an integer, and
 let $V_1,\ldots, V_m$ be a partition of
 $[n]:=\{1,\ldots,n\}$ such that $|V_k|=n\nu_k$,
 where $\nu_k$ might be the function
 of $n$ and $k=1,\ldots,m$.
 We consider the random matrix $\mathbf{A}_n(\nu_1,\ldots,\nu_m)$
 (or $\mathbf{A}_n$ for short) satisfying
 the following properties:\begin{itemize}

 \item $a_{ij}$'s are independent random variables with $a_{ij}=a_{ji}$;

 \item the $a_{ij}$'s have the same distribution $F_1$ if
 $i$ and $j\in V_k$,
  while the $a_{ij}$'s are to
 possess
 the same distribution $F_2$ if
 $i\in V_k$
 and $j\in [n]\setminus V_k$, where $k$ is an integer with $1\le k\le m$;

  \item $|a_{ij}|\le K$.
\end{itemize}
 Set
$$\mathbf{B}_n=\frac{1}{2\sqrt{n}}\mathbf{A}_n.$$ Let $\Psi_n(x)$ be the ESD of
$\mathbf{B}_n$. In section 2, we establish the LSD of $\mathbf{B}_n$
for special random matrices $\mathbf{A}_n$.
\begin{thm}\label{Main Thm}Let $F_1$ and $F_2$ be distribution functions
 mentioned above.
\begin{itemize}

 \item[\rm (i)] If
 \begin{equation}\label{Condition}\lim_{n\rightarrow\infty}\max\{\nu_1(n),\ldots,\nu_m(n)\}>0\mbox{ and }
 \lim_{n\rightarrow\infty}\frac{\nu_i(n)}{\nu_j(n)}= 1
\mbox{ for all }1\le i,j\le m,\end{equation}
  then
 $$\lim_{n\rightarrow\infty} \Psi_n(x)=\Psi(x)\mbox{ a.s.,}$$
 where $\Psi(x)$ is the limiting spectral distribution
 with density
 $$\label{SC Law-0}
 \psi(x)=\left\{\begin{array}{ll}\displaystyle
  \frac{2m}{\pi(\si_1^2+(m-1)\si_2^2)}
    \sqrt{\frac{\si_1^2+(m-1)\si_2^2}{m}-x^2} &
    \mbox{if~~ } |x|\le \sqrt{\frac{\si_1^2+(m-1)\si_2^2}{m}}_,\\
  0 & \mbox{if~~ } |x|>\sqrt{\frac{\si_1^2+(m-1)\si_2^2}{m}}_. \end{array}\right.
  $$

 \item[\rm (ii)] If
  \begin{equation}\label{con2}
  \lim_{n\rightarrow\infty}\max\{\nu_1(n),\ldots,\nu_m(n)\}=0,
  \end{equation} then
 $\lim_{n\rightarrow\infty} \Psi_n(x)=\Phi(x)\mbox{ a.s.}$
 \end{itemize}
 \end{thm}

\noindent{\bf Remark.} We require $|a_{ij}|\le K$ here for some
fixed integer $K$. In
 fact, one can readily obtain the same LSD
  for more general distributions
 $F_1$ and $F_2$ satisfying that $\si_1^2<\infty$ and $\si_2^2<\infty$
  by employing the
 classical truncation method (see \cite{B} for
 instance).\vspace{2mm}

 We then show that for general case, the random matrix $\mathbf{A}_n$ has
 no such  LSD in section 3, and, besides, propose a conjecture
 concerning the LSD of $\mathbf{A}_n$. Finally, we give an application
 about the energy of a simple graph.

\section{Proof of Theorem \ref{Main Thm}}
The main goal of this section is to show Theorem \ref{Main Thm}.
Since we can centralize the general distribution functions $F_1$ and
$F_2$, we first prove Theorem \ref{Main Thm} on condition that the
expectations $\mu_1$ and $\mu_2$ are equal to zero, and then prove
the theorem for general distributions in subsection 2.2.

\subsection{LSD for centralized distributions}

 In this part, we employ the moment approach to prove Theorem
 \ref{Main Thm} supposing that $\mu_1=\mu_2=0$.

 Above all, we deal with the first part of Theorem \ref{Main Thm}.
  It is turned out that we need to prove that the moments
 $M_{k,n}=\int_{-\infty}^{\infty}
 \la^k d\Psi_n$ $(k=1,2,\ldots)$ satisfies almost surely (a.s.)
  the following
 condition: \begin{equation}\label{Moments}
 \lim_{n\rightarrow\infty} M_{k,n}=\gamma_k=
 \left\{\begin{array}{ll} 0, & \mbox{if } k \mbox{ is odd,}\\
 \frac{k!}{2^k(k/2)!(k/2+1)!}
 f(m,\si_1,\si_2)^{k},
 & \mbox{if } k \mbox{ is even,}\end{array}\right.\end{equation}
 where
 $f(m,\si_1,\si_2)=\sqrt{\frac{\si_1^2+(m-1)\si_2^2}{m}}_.$
 Let $X_{\Psi}$ be a random variable with the distribution $\Psi(x)$. Then,
 by the linearity of expectation, we have
 $$\E(e^{itX_{\Psi}})=\sum_{k\ge 0}\frac{(it)^k}{k!}\E(X_{\Psi}^k).$$ On the
 other hand
 $$\begin{array}{lll}\displaystyle\sum_{k\ge0}\frac{(it)^k}{k!}\gamma_k
   &=&\displaystyle\sum_{j\ge0}\frac{(-1)^{j}}{j!(j+1)!}
      \left(\frac{t}{2}f(m,\si_1,\si_2)\right)^{2j}
      =\frac{2}{t\cdot f(m,\si_1,\si_2)}J_1(t\cdot f(m,\si_1,\si_2))\\
   &=&\displaystyle \frac{2}{\pi\cdot f(m,\si_1,\si_2)^2}
   \int_{-f(m,\si_1,\si_2)}^{f(m,\si_1,\si_2)}e^{itx}
   \sqrt{f(m,\si_1,\si_2)^2- x^2}~dx=\E(e^{itX_{\Psi}}),\end{array}$$
   where $J_1$ denotes the Bessel function of order 1 of the
   first kind. Therefore, $$\lim_{n\rightarrow
   \infty}M_{k,n}=\gamma_k=\E(X_{\Psi}^k)\mbox{ a.s.}, ~~k=1,2,\ldots,$$ and thus
   $\Psi_n(x)\rightarrow \Psi(x)$  a.s. $(n\rightarrow\infty)$ according to
   the Moment Convergence Theorem (\cite{B}, pp. 613).

 In order to show that $M_{k,n}\rightarrow\gamma_k$ a.s.
 $(n\rightarrow\infty)$, we first prove that
   $$\lim_{n\rightarrow\infty}\E(M_{k,n})=\gamma_k,$$ and then
   prove that
   $$\lim_{n\rightarrow\infty}\big(M_{k,n}-\E(M_{k,n})\big)=0~\mbox{a.s.}$$

   We now proceed with the calculation of $M_{k,n}$. It is not
difficult to see that
 $$\begin{array}{lll} M_{k,n} &=& \displaystyle
 \int_{-\infty}^{\infty}
 \la^k d\Psi_n=n^{-1}\sum_{j=1}^n\la_{j,n}^k=
 n^{-1}\mbox{tr}(\mathbf{B}_n^k)\\
 &=& \displaystyle2^{-k}n^{-1-k/2}\sum_{i_1=1}^n\cdots\sum_{i_k=1}^n
 a_{i_1i_2}a_{i_2i_3}\ldots a_{i_ki_1}.\end{array}$$

For $1\le v\le k$, denote by $S_{v,k,n}$ the sum of
$2^{-k}n^{-1-k/2}\E(a_{i_1i_2}a_{i_2i_3}\ldots a_{i_ki_1})$ over all
sequences $i_1,\ldots, i_k$ where $v:=\#\{i_1,\ldots, i_k\}$ (not
counting multiplicities) is the order of a sequence. Since the
expectation of $a_{ij}$ equals zero, if some $a_{ij}$ in the product
$a_{i_1i_2}a_{i_2i_3}\ldots a_{i_ki_1}$ has multiplicity one then
the expectation of the product is zero. According to the pigeon hole
principle, if $v>k/2+1$ then $S_{v,k,n}=0$, and thus
 $$\E(M_{k,n})=\sum_{v=1}^{k/2+1} S_{v,k,n}.$$
 Notice that a product
$a_{i_1i_2}a_{i_2i_3}\ldots a_{i_ki_1}$
 corresponds to an unique closed walk
 $$(i_1,i_2)(i_2,i_3)\ldots (i_k,i_1)$$ of length $k$ in the
 complete graph $K_n$ on the set $[n]$ ($K_n$ can
 contain loops here). A closed walk $(i_1,i_2)(i_2,i_3)\ldots (i_k,i_1)$
  is said to be {\em good} if $\E(a_{i_1i_2}a_{i_2i_3}\ldots a_{i_ki_1})\neq 0.$
   Hence the estimation of $S_{v,k,n}$ relies
 on a bound on the number of good walks. Let
 $W_{v,k,n}$ be the number of good walks in $K_n$ of length $k$
 and order $v$. Clearly, there are $n(n-1)\cdots(n-v+1)$ ways to
 fix an ordered good walk of $v$ distinct vertices in $K_n$.
 Moreover, for a fixed order, $W_{v,k,n}$ is a function $g(v,k)$ of
 variables $v$
 and $k$, and thus
 $$\label{Equ-1} W_{v,k,n}=n(n-1)\cdots(n-v+1)\cdot g(v,k).
 $$

 For odd $k$, since $(i_1,i_2)(i_2,i_3)\ldots (i_k,i_1)$ is a closed
 walk and $\E(a_{ij})=0$, $1\le i\le j\le n$, it is easily seen that
 $v\le (k-1)/2$ if $(i_1,i_2)(i_2,i_3)\ldots (i_k,i_1)$ is a good walk.
 Therefore, for
 $v=1,\ldots,(k-1)/2$,
 $$\begin{array}{lll} S_{v,k,n} &=& \displaystyle
 2^{-k}n^{-1-k/2}\sum_{i_1=1}^n\cdots\sum_{i_k=1}^n
 \E\left(a_{i_1i_2}a_{i_2i_3}\ldots a_{i_ki_1}\right)\\
  &\le& \displaystyle
 2^{-k}n^{-1-k/2}n^vg(v,k)K^{k}\\
 &\le& \displaystyle
 2^{-k}n^{-1}g(v,k)K^{k}\rightarrow 0~~
 (n\rightarrow\infty).
 \end{array} $$
 For even $k$, since $v\le k/2+1$,
  we again find, by a similar way, that
 $S_{v,k,n}\rightarrow 0$ as $n\rightarrow\infty$ when $v< k/2+1$.
 For the case that $v=k/2+1$, set
 $$T_{k/2}=g(k/2+1,k),$$
 {\it i.e.,} $T_k$ denotes the number of good walks $W$ in $K_n$
  of length $2k$ and
 order $k+1$ such that the order of the vertices appearing in $W$
 is fixed. We use $T'_k$ to denote the number of good walks
 $W=(i_1,i_2)(i_2,i_3)\ldots (i_{2k},i_1)$
  which contain no $i_1$ except the first and the last
 member. It is easy to see that
 $$T'_k=T_{k-1}, ~~T'_1=T_0=1, $$
 and
  \begin{equation}\label{Equ-3}
  T_k=\sum_{j=1}^k T'_jT_{k-j}=\sum_{j=1}^k T_{j-1}T_{k-j}=
  \sum_{i=0}^{k-1} T_{i}T_{k-1-i}, ~~k=1,2,\ldots. \end{equation}
  The generating function of $T_k$ is defined below
  $$T(x)=\sum_{k\ge 0}T_kx^k.$$ The recursive formula
  (\ref{Equ-3}) then gives $$T(x)=1+xT(x)^2.$$ It follows that
  $$T(x)=(2x)^{-1}\big(1\pm (1-4x)^{1/2}\big).$$ Since $T_0=1$,
  we have $T(x)=(2x)^{-1}\big(1- (1-4x)^{1/2}\big).$ Thus,
  $$T_k=\frac{1}{2}{\frac{1}{2}\choose k+1}(-4)^{k+1}
  =\frac{(2k)!}{k!(k+1)!}_{\textstyle .}$$

  To calculate $S_{k/2+1,k,n}$, we need to estimate the quantity
 that \begin{equation}\label{Quantity}\sum_{i_1=1}^n\cdots\sum_{i_k=1}^n
 \E\left(a_{i_1i_2}a_{i_2i_3}\ldots
 a_{i_ki_1}\right),\end{equation}
 where $a_{i_1i_2}a_{i_2i_3}\ldots
 a_{i_ki_1}$ corresponds to a closed walk of order $k/2+1$. In order to avoid the tedious analysis, we further assume that
  \begin{equation}\label{condition}\lim_{n\rightarrow\infty}\max\{\nu_1(n),\ldots,\nu_m(n)\}>0\mbox{ and }\nu_i=\nu_j\mbox{ for all }1\le i,j\le m.\end{equation}
  Indeed, one can readily obtain the same estimation for
  (\ref{Quantity}) on condition (\ref{Condition}) by the trick we employ below.

 Obviously, to estimate (\ref{Quantity}),
   the crucial step is to estimate the
  sum of expectations of good walks of order $k/2+1$ and
  length $k$. In this case, it is readily seen that
  each edge appears exactly twice in a good walk.
 For a fixed order of vertices appearing in a good walk
 $(i_1,i_2)(i_2,i_3)\ldots(i_k,i_1)$ of order
 $k/2+1$, we get a term $$(\si_1^2)^r(\si_2^2)^{k/2-r}T_{k/2},$$
 where $r$ is an integer with $0\le r\le k/2$.
 An edge $(i_j,i_{j+1})$ in a good walk of order $k/2+1$ and length $k$
 is said to be {\em secondary} if $i_j,i_{j+1}\in V_t$ for some part $V_t$ of the partition
 $V_1\cup\cdots\cup V_m=[n]$, otherwise, the edge is {\em chief}.
 We then pick
 up vertices according to the fixed positions of chief and
 secondary edges
  when $n$ is large enough. By the condition (\ref{condition}), we
 have the term, for large enough $n$,
 $$n(\si_1^2)^r\left(\frac{n}{m}\right)^{r}
    (\si_2^2)^{k/2-r}\left(\frac{(m-1)n}{m}\right)^{k/2-r}
    T_{k/2}
    =n^{1+k/2}T_{k/2}\left(\frac{\si_1^2}{m}\right)^r
    \left(\frac{(m-1)\si_2^2}{m}\right)^{k/2-r}.
    $$
 For any fixed value of $r$, we next choose the possible positions
 for chief and secondary edges, and thus get the term
 $$n^{1+k/2}T_{k/2}{k/2\choose r}\left(\frac{\si_1^2}{m}\right)^r
    \left(\frac{(m-1)\si_2^2}{m}\right)^{k/2-r}.$$
 According to the condition (\ref{condition}),
 $r$ may take any value from $\{0,1,\ldots,
 k/2\}$ when $n$ is large enough.
 Hence, we finally obtain the estimation of (\ref{Quantity}) that
 $$
 n^{1+k/2}T_{k/2}\sum_{r=0}^{k/2}{k/2\choose r}\left(\frac{\si_1^2}{m}\right)^r
    \left(\frac{(m-1)\si_2^2}{m}\right)^{k/2-r}
   =n^{1+k/2}T_{k/2}\left(\frac{\si_1^2}{m}+\frac{(m-1)\si_2^2}{m}\right)^{k/2}.$$
Therefore, for large enough $n$, we have $$\begin{array}{lll}
S_{k/2+1,k,n} &=& \displaystyle
 2^{-k}n^{-1-k/2}\sum_{i_1=1}^n\cdots\sum_{i_k=1}^n
 \E\left(a_{i_1i_2}a_{i_2i_3}\ldots a_{i_ki_1}\right)\\
 &=& \displaystyle 2^{-k}n^{-1-k/2}\cdot
   n^{1+k/2}T_{k/2}\left(\frac{\si_1^2+(m-1)\si_2^2}{m}\right)^{k/2}\\
 &=& \displaystyle 2^{-k}\frac{k!}{(k/2)!(k/2+1)!}
 \left(\frac{\si_1^2+(m-1)\si_2^2}{m}\right)^{k/2}\\
 &=& \displaystyle
 \frac{k!}{2^k(k/2)!(k/2+1)!}f(m,\si_1,\si_2)^k=\gamma_k.
 \end{array} $$ Consequently,
  $$\E(M_{k,n})=S_{k/2+1,k,n}\rightarrow \gamma_k ~~~(n\rightarrow\infty).$$

 We now estimate the difference $M_{k,n}-\E(M_{k,n})$. Using
 Markov's inequality, we have
 $$\label{Difference }
 \Pro[|M_{k,n}-\E(M_{k,n})|>\ep]\le \E\big[(M_{k,n}-\E(M_{k,n}))^2\big]/\ep^2.
 $$
 Hence, to prove that
$\Pro\left[\lim_{n\rightarrow\infty}\big(M_{k,n}-\E(M_{k,n})\big)=0\right]=1$,
it suffices to show that
\begin{equation}\label{The series}\sum_{n\ge1}\E\big[(M_{k,n}-\E(M_{k,n}))^2\big]<\infty,
\mbox{ for any given }k.\end{equation}
 One can readily see that
 $$\begin{array}{ll}
 &\E\big[(M_{k,n}-\E(M_{k,n}))^2\big]\\
  =& \E\big[(M_{k,n})^2\big]-(\E[M_{k,n}])^2\\
 =&\displaystyle 2^{-2k}n^{-2-k}\sum_{i_1=1}^n\cdots\sum_{i_k=1}^n
 \sum_{j_1=1}^n\cdots\sum_{j_k=1}^n \\
 &\Big[
 \E(a_{i_1i_2}\ldots a_{i_ki_1}a_{j_1j_2}\ldots a_{j_kj_1})-
 \E(a_{i_1i_2}\ldots a_{i_ki_1})\E(a_{j_1j_2}\ldots
 a_{j_kj_1})\Big].
 \end{array}$$
 We use $a_i$ and $a_j$ to denote,
respectively, the sequences $a_{i_1i_2}\ldots a_{i_ki_1}$ and
$a_{j_1j_2}\ldots a_{j_kj_1}$. Obviously, to prove (\ref{The
series}), it suffices to show that if $\E(a_i\cdot a_j)-
 \E(a_i)\E(a_j)\neq 0
 $, then $|V(W_i)\cup V(W_j)|\le k$, where $W_i$ and $W_j$
 denote the two closed walks $W_i:=(i_1,i_2)\ldots (i_k,i_1)$ and
 $W_j:=(j_1,j_2)\ldots(j_k,j_1)$, respectively.

 One can easily see that if $a_i$ and $a_j$ are independent or
 $\E(a_i\cdot a_j)=0$ then \linebreak $\E(a_i\cdot a_j)-
 \E(a_i)\E(a_j)=0$. Thus it is sufficient to consider the case
 that $a_i$ and $a_j$ are not independent and
 $\E(a_i\cdot a_j)\neq0.$
 \begin{claim} If $a_i$ and $a_j$ are not independent and
 $\E(a_i\cdot a_j)\neq0,$ then $|V(W_i)\cup V(W_j)|\le k$.
\end{claim}
  Clearly, if $a_i$ and $a_j$ are not independent
 then $V(W_i)\cap V(W_j)\neq\emptyset$. Then $W_i\cup W_j$ is a
 closed walk of length $2k$ since $W_i$ and $W_j$ are two closed walks of
 length $k$, respectively. If
 $\E(a_i\cdot a_j)\neq0$
 then the order of $W_i\cup W_j$ is not more than $k+1$ by pigeon
 hole principle. Furthermore, if $|V(W_i)\cup V(W_j)|=k+1$ then
 $a_i$ and $a_j$ are independent. In fact, each edge in $W_i\cup
 W_j$ appears exactly twice when $|V(W_i)\cup V(W_j)|=k+1$.
 Thus, those edges induce (not counting multiplicities ) a tree in
 $K_n$ of order $k+1$ since $W_i\cup W_j$ is connected. We
 further assert that $E(W_i)\cap E(W_j)=\emptyset.$ Suppose, for
 a contradiction, that there exists one element
 $a_{i_s}a_{i_{s+1}}$ appears only once in $W_i$. Then the
 subgraph graph induced by $W_i$ should contain a cycle since
 $W_i$ is a closed walk, which contradicts to the fact that the
 subgraph induced by $W_i\cup W_j$ is a tree. Therefore, $E(W_i)\cap
 E(W_j)$ is empty, and thus $a_i$ and $a_j$ are independent.
 Hence, our claim follows.

 We thus have $$\E\big[(M_{k,n}-\E(M_{k,n}))^2\big]\le n^{-2},$$
 and then (\ref{The series}) holds. Therefore, (\ref{Moments}) follows, and
 this completes our proof of the
 first part of
 Theorem \ref{Main Thm} on condition that $\mu_1=\mu_2=0$.

 We next show the second part of Theorem \ref{Main Thm}  on condition that $\mu_1=\mu_2=0$.
  We can hold the desire by applying the moment approach again. In fact,
  by the approach, it is sufficient to show that
 the moments
 $M_{k,n}$ $(k=1,2,\ldots)$  satisfies a.s.
  the following
 condition: \begin{equation}\label{moments}
 \lim_{n\rightarrow\infty} M_{k,n}=\gamma_k=
 \left\{\begin{array}{ll} 0, & \mbox{if } k \mbox{ is odd,}\\
 \frac{k!}{2^k(k/2)!(k/2+1)!}
 \si_2^{k},
 & \mbox{if } k \mbox{ is even.}\end{array}\right.\end{equation}
 It is
 similar to the proof of the first part that the crucial step
 is to estimate the quantity (\ref{Quantity}) when $k$ is even.
 Evidently, in this case, each edge appears exactly twice in a good walk.
 Let $\mathcal{W}'_{k/2+1,k,n}$ be the set of good walks
 in $K_n$ of order $k/2+1$
 and length $k$ in which each walk contains at least one secondary edge
 (not counting multiplicities ). Set
 $$W'_{k/2+1,k,n}=\big|\mathcal{W}'_{k/2+1,k,n}\big|\mbox{ and }
 W''_{k/2+1,k,n}=W_{k/2+1,k,n}-W'_{k/2+1,k,n}.$$
  It is not hard to see that
 $$
   W'_{k/2+1,k,n}
      \le \frac{k}{2}\cdot n^{k/2} o(n) T_{k/2}
   =o(n^{k/2+1}).
 $$
 Thus,
 $$
 \lim_{n\rightarrow\infty}\frac{W_{k/2+1,k,n}}{n^{k/2+1}}
 =\lim_{n\rightarrow\infty}\frac{W'_{k/2+1,k,n}+W''_{k/2+1,k,n}}{n^{k/2+1}}
 =\lim_{n\rightarrow\infty}\frac{W''_{k/2+1,k,n}}{n^{k/2+1}}_{\displaystyle ,}$$
 and $$\sum_{W_i\in\mathcal{W}'_{k/2+1,k,n}}\E(W_i)
 \le o(n^{k/2+1})K^k.$$
 On the other hand,
$$W_{k/2+1,k,n}=n(n-1)\cdots(n-k/2)\cdot T_{k/2}
=n(n-1)\cdots(n-k/2)\frac{(2k)!}{k!(k+1)!}_. $$
 Hence, for large enough $n$, we have $$\begin{array}{lll}
S_{k/2+1,k,n} &=& \displaystyle
 2^{-k}n^{-1-k/2}\sum_{i_1=1}^n\cdots\sum_{i_k=1}^n
 \E\left(a_{i_1i_2}a_{i_2i_3}\ldots a_{i_ki_1}\right)\\
 &=& \displaystyle 2^{-k}n^{-1-k/2}\cdot
   W''_{k/2+1,k,n}(\si_2^2)^{k/2}\\
 &=& \displaystyle 2^{-k}n^{-1-k/2}\cdot
   W_{k/2+1,k,n}(\si_2^2)^{k/2}\\
 &=& \displaystyle\frac{k!}{2^k(k/2)!(k/2+1)!}
 \si_2^{k}=\gamma_k.
 \end{array} $$
 Therefore, $$\lim_{n\rightarrow\infty}\E(M_{k,n})=\gamma_k.$$
 Using a similar way in the proof of the first part, one can also
 prove that $$M_{k,n}\rightarrow\E(M_{k,n})
 ~a.s.~(n\rightarrow\infty).$$ Thus (\ref{moments}) holds and the second part of
 Theorem \ref{Main Thm} follows when
 $\mu_1=\mu_2=0$.

\subsection{LSD for general distributions}

In this subsection, we show that Theorem \ref{Main Thm} holds for
general distribution functions $F_1$ and $F_2$ by two distinct
tools.

In the following, the norm $||f||$ of a real function $f$ is always
defined as follows:
$$||f||=\sup_x|f(x)|.$$
\begin{lem}[\bf Rank Inequality \cite{B}]\label{Rank Ineq}
Let $\mathbf{U}$ and $\mathbf{V}$ be two Hermitian matrices of order
$n$, and let $\Psi_\mathbf{U}(x)$ and $\Psi_\mathbf{V}(x)$ be the
ESD of $\mathbf{U}$ and $\mathbf{V}$, respectively. Then
$$\left|\left|\Psi_\mathbf{U}(x)-\Psi_\mathbf{V}(x)\right|\right|\le\frac{1}{n}\rank(\mathbf{U}-\mathbf{V}).$$
\end{lem}
The {\em Stieltjes transform} $S(z)$ of a function $F(x)$ is defined
below
$$S(z)=\int_{-\infty}^{\infty}(x-z)^{-1}~dF(x),\ \textmd{Im}(z)>0.$$
One can readily see that for the ESD $\Psi_n(x)$ of $\mathbf{B}_n$,
we have
$$\int_{-\infty}^{\infty}(x-z)^{-1}~d\Psi_n(x)=n^{-1}\textmd{tr}(\mathbf{B}_n-zI)^{-1}.$$
Here, we need two facts about this transform, and refer the readers
to \cite{stj} for details.

\begin{lem}\label{cst}~
\begin{itemize}

\item[{\rm (i)}] $F(x)$ is uniquely determined by $S(z)$.

\item[{\rm (ii)}] For probability distribution, $F_n(x)\rightarrow F(x)$ if
and only if $S_n(z)\rightarrow S(z)$ pointwise.

\end{itemize}
\end{lem}

Let $\mathbf{A}_n$ be a symmetric matrix, and let $\mathbf{D}_n$ be
a symmetric quasi-diagonal matrix. We use $\overline{\Psi}_n(x)$ to
denote the ESD of $\mathbf{A}_n+\mathbf{D}_n$. Then we have the
following result.
\begin{lem}\label{st}
Let $S_n(z)$ and $\overline{S}_n(z)$ be the Stieltjes transforms of
$\Psi_n(x)$ and $\overline{\Psi}_n(x)$, respectively. Then
$$|S_n(z)-\overline{S}_n(z)|\leq {\mbox{\rm Im}}(z)^{-2}||\mathbf{D}_n||_1,
$$
where $||\mathbf{D}_n||_1:=\max_{j\in[n]}\{\sum_{i=1}^nD_n(ij)\}$ is
the 1-normal number of $\mathbf{D}_n$.
\end{lem}

Denote by $\lambda(\mathbf{M})$ the spectral radius for some real
symmetric matrix $\mathbf{M}$ of order $n$. Clearly,
$\frac{1}{n}|\mbox{tr}(\mathbf{M})|\leq
\lambda(\mathbf{M})\leq||\mathbf{M}||_1$. As is well known, the
eigenvalues of $\mathbf{M}$ are real. Then
$\lambda(\mathbf{M}-zI)\geq |\mbox{Im}(z)|.$ By these observations,
we show lemma \ref{st} as follows.

\begin{proof}[\bf Proof] Clearly,
$(\mathbf{A}_n-zI)^{-1}-(\mathbf{A}_n+\mathbf{D}_n-zI)^{-1}=(\mathbf{A}_n+\mathbf{D}_n-zI)^{-1}\mathbf{D}_n(\mathbf{A}_n-zI)^{-1}.$
 Then
\begin{eqnarray*}
|S_n(z)-\overline{S}_n(z)|&\leq&n^{-1}\Big|\mbox{tr}\big((\mathbf{A}_n+\mathbf{D}_n-zI)^{-1}\mathbf{D}_n(\mathbf{A}_n-zI)^{-1}\big)\Big|\\
&\leq&\lambda\big((\mathbf{A}_n+\mathbf{D}_n-zI)^{-1}\mathbf{D}_n(\mathbf{A}_n-zI)^{-1}\big)\\
&\leq&(\mbox{Im}(z))^{-2}||\mathbf{D}_n||_1.
\end{eqnarray*}\end{proof}

We assume, without loss of generality, that $\mathbf{A}_n$ is a
random matrix with the partition $V_1,\ldots,V_m$ such that
$n\nu_i\rightarrow\infty$, $i=1,2,\ldots,l$, and $n\nu_i<\infty$,
$i=l+1,\ldots,m$, as $n\rightarrow\infty$. Let $\mathbf{H}_n$ be a
quasi-diagonal matrix of order $n$ such that
$$h_{ij}=\left\{\begin{array}{ll}
 1, & \mbox{if } i \mbox{ and } j\in V_k(1\le k\le m),\\
 0, & \mbox{otherwise,}
\end{array}\right.$$
and let $\mathbf{H}_n'$ be a matrix such that
$$h'_{ij}=\left\{\begin{array}{ll}
 1, & \mbox{if } i \mbox{ and } j\in V_k (k\le l),\\
 0, & \mbox{otherwise. }
\end{array}\right.$$
Set
$$\mathbf{H}_n''=\mathbf{H}_n-\mathbf{H}_n',$$
$$\mathbf{C}_n'=\frac{1}{2\sqrt n}(\mathbf{A}_n-(\mu_1-\mu_2)\mathbf{H}_n'-\mu_2\mathbf{J}_n),$$
and
$$\mathbf{C}_n''=\mathbf{C}_n'-\frac{1}{2\sqrt n}((\mu_1-\mu_2)\mathbf{H}_n''),$$
where $\mathbf{J}_n$ is the  matrix in which all elements equal 1.

Since $\E(\mathbf{C}_n''(ij))=0$, the LSD of $\mathbf{C}_n''$ is
$\Psi(x)$ if (\ref{Condition}) holds (or $\Phi(x)$ if (\ref{con2})
holds), as shown above. Let $\Psi_n''$ and $\Psi_n'$ be the ESD of
$\mathbf{C}_n''$ and $\mathbf{C}_n'$, respectively. Then, for their
corresponding Stieljes transforms $S_n''(z)$ and $S_n'(z)$, we have
from Lemma \ref{st} that
$$|S_n''(z)-S_n'(z)|\leq (\mbox{Im}(z))^{-2}\frac{1}{2\sqrt
n}||(\mu_1-\mu_2)\mathbf{H}_n''||_1.
$$
Since each block matrix on the diagonal of $\mathbf{H}_n''$ is of
finite order, we have $$\frac{1}{2\sqrt
n}||(\mu_1-\mu_2)\mathbf{H}_n''||_1\rightarrow 0
~(n\rightarrow\infty).$$ Then, we can get that
$\lim_{n\rightarrow\infty}S_n''(z)=\lim_{n\rightarrow\infty}S_n'(z)$
for any $z$ such that $\mbox{Im}z>0$. Because $\Psi(x)$ (or
$\Phi(x)$) is the LSD of $\mathbf{C}_n''$, from Lemma \ref{cst}(ii),
we have that $\lim_{n\rightarrow\infty}S_n''(z)=S^*(z)$, where
$S^*(z)$ is the Stieltjes transform of $\Psi(x)$ (or $\Phi(x)$).
Therefore, $\lim_{n\rightarrow\infty}S_n'(z)=S^*(z)$, and thus the
ESD $\Psi_n'(x)$ of $\mathbf{C}_n'$ converges to $\Psi(x)$ (or
$\Phi(x)$). So, $\mathbf{C}_n'$ has the same LSD as
$\mathbf{C}_n''$.

Furthermore,  since each of the block matrices on the diagonal of
$\mathbf{H}_n'$ is of infinite order, there are $o(n)$ such block
matrices. Then we have
$\mbox{rank}|(\mu_1-\mu_2)\mathbf{H}_n'+\mu_2\mathbf{J}_n|=o(n)$. By
employing Lemma \ref{Rank Ineq} for $\mathbf{B}_n$ and
$\mathbf{C}_n'$, $\Psi_n(x)$ of $\mathbf{B}_n$ converges to the LSD
$\Psi(x)$ (or $\Phi(x)$).

Therefore, Theorem \ref{Main Thm} holds for general distribution
functions $F_1$ and $F_2$.

\section{The LSD for more general random matrices}

In this section, we shall show that there is no LSD for general
random matrix $\mathbf{A}_n(\nu_1,\ldots,\nu_m)$. Actually, we shall
prove that $\mathbf{A}_n(\nu_1,\ldots,\nu_m)$ has no LSD for some
special cases. To be precise, we shall show that if $F_1\equiv 0$
then the LSD $\mathbf{A}_n(\nu_1,\nu_2)$ exists if and only if
$\lim_{n\rightarrow\infty}\nu_1/\nu_2=1$. Moreover, we verify that
there exist real numbers $\nu_1,\ldots,\nu_m$ such that
$\mathbf{A}_n(\nu_1,\ldots,\nu_m)$ has no LSD under the condition
$F_1\equiv 0$. In general, we have the following conjecture.
\begin{con}\label{conj}
Let $\mathbf{A}_n$ be a random matrix with partition
$V_1,\ldots,V_m$ such that\linebreak
$\lim_{n\rightarrow\infty}\max\{\nu_1(n),\ldots,\nu_m(n)\}>0$. If
the LSD of $\mathbf{A}_n$ exists, then
$\lim_{n\rightarrow\infty}\frac{\nu_i(n)}{\nu_j(n)}= 1$ for all
$i,j\in[m]$.
\end{con}
\noindent{\bf Remark.} According to Theorem \ref{Main Thm}, if
$\lim_{n\rightarrow\infty}\nu_i(n)/\nu_j(n)= 1$ for all $i,j\in[m]$
then the LSD of $\mathbf{A}_n$ exists. Thus, if Conjecture
\ref{conj} is true, then the condition that
$\lim_{n\rightarrow\infty}\frac{\nu_i(n)}{\nu_j(n)}= 1$ for all
$i,j\in[m]$ is necessary and sufficient for the existence of the LSD
of $\mathbf{A}_n$.

Firstly, we investigate the LSD for $\mathbf{A}_n(\nu_1,\nu_2)$ when
$F_1\equiv0$ and $0<\nu_1<1$.  A function $f(t)$ is defined to be
{\em nonnegative} if it satisfies that
$$\sum_{k=1}^l\sum_{j=1}^lf(t_k-t_j)r_k\overline{r}_j\geq0,$$
for any positive integer $l$, and any real numbers $t_1,\ldots,t_l$
and complex numbers $r_1,\ldots,r_l.$ In order to avoid the tedious
analysis, we further assume that $\nu_1$ and $\nu_2$ are real
numbers. Indeed, one can readily obtain the same result under the
condition $\lim_{n\rightarrow\infty}\nu_1/\nu_2=1$ by the method we
employ below.
\begin{thm}\label{p1}
Let $\mathbf{A}_n(\nu_1,\nu_2)$ be a random matrix with $F_1\equiv0$
and $0<\nu_1<1$. If $\nu_1\neq\nu_2$, then
$\mathbf{A}_n(\nu_1,\nu_2)$ has no LSD.
\end{thm}
\begin{proof}[\bf Proof]
We prove the assertion  by a contradiction. Since we can centralize
the general distribution $F_2$, we further assume that
$\E(a_{ij})=0$ in what follows. We suppose that $\nu_1\neq\nu_2$ and
there exists a function $\Psi(x)$ such that the ESD $\Psi_n(x)$ of
$\mathbf{B}_n=\mathbf{A}_n(\nu_1,\nu_2)/(2\sqrt{n})$ converges to it
almost surely as $n\rightarrow\infty$. Then
$M_{k,n}=\int_{-\infty}^{\infty}x^k~d\Psi_n$ converges to
$\gamma_k=\int_{-\infty}^{\infty}x^k~d\Psi$ almost surely as
$n\rightarrow\infty$ $(k=1,2,\ldots)$. So we can get the estimation
of $\gamma_k$ by calculating the moment $M_{k,n}$.

We first estimate $\E(M_{k,n})=\sum_{v=1}^{k/2+1} S_{v,k,n}$, as we
did in subsection 2.1. For the similar reason, one can readily see
that we merely need to compute $\E(M_{k,n})$ for even $k=2j$.
Moreover, if $v< j+1$ then $S_{v,2j,n}\rightarrow 0$ as
$n\rightarrow\infty$. Thus, to get the estimation of $\E(M_{k,n})$,
it suffices to focus on $S_{j+1,2j,n}$ which indeed satisfies that
$$\lim_{n\rightarrow \infty}S_{j+1,2j,n} = \left\{\begin{array}{ll}
 2^{-2j}\cdot2(\nu_1\nu_2)^{j/2}\sigma_2^{2j}\cdot
T_j & \mbox { if  $j$ is even},\\
2^{-2j}\cdot(\nu_1\nu_2)^{(j-1)/2}\sigma_2^{2j}\cdot T_j & \mbox {
if  $j$ is odd}.
\end{array}\right.$$
One then can prove by an analogous way in Subsection 2.1 that
$$\lim_{n\rightarrow\infty}\big(M_{k,n}-\E(M_{k,n})\big)=0\mbox{ a.s.}$$
Set $\widehat{\nu}=(\nu_1\nu_2)^{1/4}.$ Apparently,
$0<\widehat{\nu}<\sqrt{1/2}$ since $\nu_1\neq\nu_2$. Then, the
following equality holds a.s.
\begin{eqnarray*}
\label{liy1}
\lim_{n\rightarrow\infty}M_{k,n}=\left\{\begin{array}{lll}
0  &k\equiv 1,3\mod4, \\[8pt]
\displaystyle\frac{2k!{\widehat{\nu}}^k\sigma_2^{k}}{2^k(k/2)!(k/2+1)!} &k\equiv 0\mod4,\\[10pt]
\displaystyle\frac{k!{\widehat{\nu}}^k\sigma_2^{k}}{{\widehat{\nu}}^22^k(k/2)!(k/2+1)!}
&k\equiv 2\mod4.
\end{array}
\right.
\end{eqnarray*}
Let $X_{\Psi}$ be a random variable with the distribution $\Psi$,
and let $f(t):=\E(e^{itX_{\Psi}}) =
\sum_{k=0}^{\infty}\frac{\gamma_k}{k!}(it)^k$ be the character
function of $X_{\Psi}$.  Since $M_{k,n}\rightarrow\gamma_k$ a.s.
$(n\rightarrow\infty)$, $k=1,2,\ldots$, the following equality holds
a.s.
\begin{eqnarray*}
\E(e^{itX_{\Psi}}) &=& \sum_{\scriptsize\mbox{\it j is
even}}^{\infty}\frac{2\cdot
(-1)^j}{j!(j+1)!}\left(\frac{t}{2}\right)^{2j}{\widehat{\nu}}^{2j}\sigma_1^{2j}+\sum_{\scriptsize\mbox{\it
j is odd}}^{\infty}\frac{
(-1)^j}{{\widehat{\nu}}^2j!(j+1)!}\left(\frac{t}{2}\right)^{2j}{\widehat{\nu}}^{2j}\sigma_1^{2j}\\
&=&2\cdot
\frac{J_1({\widehat{\nu}}\sigma_Ft)+I_1(\widehat{\nu}\sigma_Ft)}{\widehat{\nu}\sigma_Ft}+\frac{1}{{\widehat{\nu}}^2}\cdot
\frac{J_1({\widehat{\nu}}\sigma_Ft)-I_1(\widehat{\nu}\sigma_Ft)}{\widehat{\nu}\sigma_Ft}\\
&=&
(2+\frac{1}{\widehat{\nu}^2})\cdot\frac{J_1(\widehat{\nu}\sigma_Ft)}{\widehat{\nu}\sigma_Ft}+(2-\frac{1}{\widehat{\nu}^2})\cdot\frac{I_1(\widehat{\nu}\sigma_Ft)}{\widehat{\nu}\sigma_Ft}\\
&=&
(2+\frac{1}{\widehat{\nu}^2})\cdot\frac{1}{\pi}\int_{-1}^{1}e^{i\widehat{\nu}\sigma_Ftx}\sqrt{1-x^2}dx+(2-\frac{1}{\widehat{\nu}^2})\cdot\frac{1}{\pi}\int_{-1}^{1}e^{-\widehat{\nu}\sigma_Ftx}\sqrt{1-x^2}dx
\end{eqnarray*}
($J_1$ denotes the Bessel function of order 1 of the first kind and
$I_1$ denotes the modified Bessel function of order 1 of the first
kind, such that $I_1(t)=-iJ_1(it)$). We further assume that $l=2,
r_1=1,r_2=1$ and $t_2=0$. Since
$0<\widehat{\nu}<\sqrt{\frac{1}{2}}$, if $t_1$ is large enough then
\begin{eqnarray*}&
&~~~\sum_{k=1}^l\sum_{j=1}^lf(t_k-t_j)r_k\overline{r}_j
\\&
&=\int_{-1}^{1}\sum_{k=1}^{2}\sum_{k=1}^{2}\left[(2+\frac{1}{\widehat{\nu}^2})e^{i\widehat{\nu}\sigma_2(t_k-t_j)x}\cdot\frac{1}{\pi}\sqrt{1-x^2}+
(2-\frac{1}{\widehat{\nu}^2})e^{-\widehat{\nu}\sigma_2(t_k-t_j)x}\cdot\frac{1}{\pi}\sqrt{1-x^2}\right]dx\\
&
&\leq\int_{-1}^{1}(2+\frac{1}{\widehat{\nu}^2})\cdot\frac{1}{\pi}\sqrt{1-x^2}\left|\sum_{k=1}^2e^{i\widehat{\nu}\sigma_2
t_kx}\right|^2~dx+\int_{-1}^{1}(2-\frac{1}{\widehat{\nu}^2})e^{\widehat{\nu}\sigma_2
t_1x}\cdot\frac{1}{\pi}\sqrt{1-x^2}~dx\\&
&\leq2(2+\frac{1}{\widehat{\nu}^2})+(2-\frac{1}{\widehat{\nu}^2})\int_{1/4}^{1/2}\frac{1}{\pi}\cdot
e^{\widehat{\nu}\sigma_2t_1/4}\cdot\frac{\sqrt{3}}{2}~dx\\
&
&\leq2(2+\frac{1}{\widehat{\nu}^2})+(2-\frac{1}{\widehat{\nu}^2})\frac{1}{\pi}\cdot
e^{\widehat{\nu}\sigma_2t_1/4}\cdot\frac{\sqrt{3}}{8}<0,
\end{eqnarray*}
which contradicts to the fact that the character function should be
nonnegative. Hence, the necessity follows.\end{proof}

For general random matrix $\mathbf{A}_n(\nu_1,\ldots,\nu_m)$,
$m\geq3$, we fail to obtain the result in the same way as the case
for $m=2$. Indeed, if we estimate the moment $\gamma_k$ by the same
method we used above, the step to count the number of good walks for
$S_{j+1,2j,n}$ is  much complicated. Worse still, the character
function is also harder to get.  However, we can still verify
Conjecture \ref{conj} for some special cases.

\begin{pro}
Suppose $F_1\equiv 0$, $\sigma_2=1$, $\nu_1>3/4$ and
$\nu_2=\cdots=\nu_m$, where $m\geq 3$. Then, there exists no LSD for
$\mathbf{A}_n(\nu_1,\ldots,\nu_m)$ a.s.
\end{pro}

\begin{proof}[\bf Proof]
Since we can centralize the general distribution $F_2$, we further
assume $\mu_2=0$ in what follows. For a contradiction, we assume
$\Psi(x)$ is the LSD of $\mathbf{B}_n$ such that
$\lim_{n\rightarrow\infty}\Psi_n(x)= \Psi(x)$ a.s. Then, $M_{k,n}=
\int x^k~d\Psi_n(x)\rightarrow \int x^k~d\Psi(x)=\gamma_k$ a.s.
$(n\rightarrow\infty)$.

It is readily seen that for any $k+1$ real numbers
$t_0,t_1,\ldots,t_k$,
\begin{equation}\label{equ-pq}
\sum_{p,q=0}^k\gamma_{p+q}t_pt_q=\int_{-\infty}^{\infty}(t_0+t_1x+\cdots+t_kx^k)^2~d\Psi(x)\geq0.
\end{equation}
Obviously, $\gamma_0=1$. Set $$\Delta_k=\left(\begin{array}{llll}
\gamma_0 & \gamma_1 & \cdots &\gamma_k \\
\gamma_1 & \gamma_2 & \cdots&\gamma_{k+1}\\
\vdots & \vdots & \cdots&\vdots\\
\gamma_k & \gamma_{k+1} & \cdots&\gamma_{2k}\\
\end{array}\right).$$
Then due to (\ref{equ-pq}), we have
$$(t_0,\cdots,t_k)\Delta_k(t_0,\cdots,t_k)^T\geq0.$$
Thus, the symmetric matrix $\Delta_k$ is non-negative definite for
any $k$. Therefore, $|\Delta_k|\geq 0$, $k=0,1,2,\cdots.$

Then, applying the same method as in Section 2.1, we can compute the
moment $\gamma_k$. At first, to estimate $\E(M_{k,n})$, we just need
to focus on the case for $k=2j$, and, moreover, it suffices to
calculate $S_{j+1,2j,n}$.

Let $j=1$. We have $$S_{2,2,n}=2^{-2}\cdot
n^{-2}\sum_{i_1=1}^n\sum_{i_2=1}^n\E(a_{i_1i_2}^2).$$ If $i_1,i_2$
are in the same part $V_i \ (1\leq i\leq m)$, $\E(a_{i_1i_2}^2)=0$
since $F_1\equiv0$. Thus, only for $i_i,i_2$ that do not belong to
the same part $V_i$, the expectation $\E(a_{i_1i_2}^2)$ contributes
1 to the value $S_{2,2,n}$.  Combining the fact
$(m-1)\nu_2+\nu_1=1$, we have $$S_{2,2,n}\rightarrow
2^{-2}\cdot(\nu_1(m-1)\nu_2+(m-1)\nu_2(1-\nu_2))~\mbox{as }
n\rightarrow\infty.$$ One then can prove by an analogous way in
Subsection 2.1 that
$$\lim_{n\rightarrow\infty}\big(M_{k,n}-\E(M_{k,n})\big)=0\mbox{ a.s.}$$
Therefore, the value
$\gamma_2=2^{-2}\cdot(\nu_1(m-1)\nu_2+(m-1)\nu_2(1-\nu_2))$ a.s.

By this means, we have
$$\begin{array}{lll}

 \gamma_2 &=& 2^{-2}\cdot(\nu_1(m-1)\nu_2+(m-1)\nu_2(1-\nu_2)),\\
 \gamma_4 &=& 2\cdot2^{-4}\cdot(\nu_1(m-1)\nu_2(1-\nu_2)+(m-1)\nu_2\nu_1(m-1)\nu_2\\
   & &+(m-1)\nu_2(m-2)\nu_2(1-\nu_2)),\\
 \gamma_6 &=& 5\cdot2^{-6}\cdot(\nu_1(m-1)\nu_2\nu_1(m-1)\nu_2+\nu_1(m-1)\nu_2(m-2)\nu_2(1-\nu_2)\\
   & &+(m-1)\nu_2\nu_1(m-1)\nu_2(1-\nu_2)+(m-1)\nu_2(m-2)\nu_2\nu_1(m-1)\nu_2\\
   & &+(m-1)\nu_2(m-2)\nu_2(m-2)\nu_2(1-\nu_2)). \end{array}$$
and $\gamma_1, \gamma_3,\gamma_5=0$.

But then $|\Delta_3|<0$, which contradicts to the fact
$|\Delta_k|\geq0$. Thus, the proposition follows.
\end{proof}

\section{Application to the energy of random graphs}

In this section, we shall compute the energy of a random graph by
the results established in the previous sections. Our notions and
terminology are standard, and we refer the readers to \cite{bb} for
the conceptions not defined here. Let $G$ be a simple graph of order
$n$. The eigenvalues $\lambda_1,\ldots,\lambda_n$ of the adjacent
matrix of $G$ are said to be the {\em eigenvalues} of $G$. In
chemistry, there is closed correspondence between the graph
eigenvalues and the molecular orbital energy levels of
$\pi$-electrons in conjugated hydrocarbons. For the H\"{u}chkel
molecular orbital (HMO) approximation, the total $\pi$-electron
energy $\En(G)$ in conjugated hydrocarbons is given by the sum of
absolute values of the eigenvalues corresponding to the molecular
graph $G$. In 1970s, Gutman \cite{gut} extended the conception of
energy to all simple graphs who defined
$$\En(G)=\sum_{i=0}^n|\lambda_i|,$$
where $\lambda_1,\ldots,\lambda_n$ are the eigenvalues of $G$.
Recently, this graph invariant has attracted a lot of attention, and
the readers are refereed to \cite{glz} for further details.

Let $G_n(p)$ be a random graph of $\mathcal{G}_{n}(p)$. It is easy
to see that if $F_1\equiv 0$ and $F_2$ is a Bernoulli distribution
with mean $p$, then the random matrix $\mathbf{X}_n$ is the adjacent
matrix of $G_n(p)$. According to Theorem \ref{Thm-1}, almost every
(a.e.) random graph $G_n(p)$ enjoys the equation below:
$$\begin{array}{lll}
  \En(G_{n}(p))
  &=&\displaystyle
  2\sqrt{n}\cdot n\left(\frac{2}{\pi\si_2^2}\int_{-\si_2}^{\si_2}
  |x|\sqrt{\si_2^2-x^2}~dx+o(1) \right)\\
  &=&\displaystyle
  n^{3/2}\left(\frac{8}{3\pi}\si_2+o(1)\right)
  =n^{3/2}\left(\frac{8}{3\pi}\sqrt{p(1-p)}+o(1)\right).
 \end{array}$$

Note that for $p=\frac 1 2$, Nikiforov in \cite{N} got the above
formula. Here, our result is for any probability $p$, which could be
seen as a generalization of his result. Next we will get the energy
for random $m$-partite graphs.

We use $K_{n;\nu_1,\ldots,\nu_m}$ to denote the complete $m$-partite
graph of order $n$ whose parts $V_1,\ldots,V_m$ are such that
$|V_i|=n\nu_i$, $i=1,\ldots,m$, where $m=m(n)\ge 2$ is an integer.
Let $\mathcal{G}_{n;\nu_1\ldots\nu_m}(p)$ be the set of random
graphs in which the edges are chosen independently with probability
$p$ from $K_{n;\nu_1,\ldots,\nu_m}$. Especially, we denote by
$K[n;m]$ and $\mathcal{G}_{n,m}(p)$, respectively, the complete
$m$-partite graph and the set of $m$-partite random graphs
satisfying
$$\lim_{n\rightarrow\infty}\max\{\nu_1(n),\ldots,\nu_m(n)\}>0\mbox{ and }
 \lim_{n\rightarrow\infty}\frac{\nu_i(n)}{\nu_j(n)}=1.
 $$
One can readily see that if a random matrix $\mathbf{A}_n$ and the
complete $m$-partite graph $K[n;m]$ have the same partition, and
$F_1\equiv 0$ and $F_2$ is a Bernoulli distribution with mean $p$,
then $\mathbf{A}_n$ is the adjacent matrix of
$G_{n,m}(p)\in\mathcal{G}_{n,m}(p)$. Employing the first part of
Theorem \ref{Main Thm}, a.e. random graph $G_{n,m}(p)$ enjoys the
following equation:
$$\begin{array}{lll}
  \En(G_{n,m}(p))
  &=&\displaystyle
  2\sqrt{n}\cdot n\left( \frac{2m}{\pi(m-1)\si_2^2}
  \int^{\sqrt{\frac{m-1}{m}}\si_2}_{-\sqrt{\frac{m-1}{m}}\si_2}
    |x|\sqrt{\frac{(m-1)\si_2^2}{m}-x^2}~dx+o(1)\right)\\
  &=& \displaystyle
  n^{3/2}\left(\frac{8}{3\pi}\sqrt{\frac{m-1}{m}}\si_2+o(1)\right)\\
  &=& \displaystyle
  n^{3/2}\left(\frac{8}{3\pi}\sqrt{\frac{m-1}{m}p(1-p)}+o(1)\right).
\end{array}$$

Furthermore, we can get the energy $\En$ of a random graph
$G_{n;\nu_1\ldots\nu_m}(p)\in\mathcal{G}_{n;\nu_1\ldots\nu_m}(p)$ if
$\lim_{n\rightarrow\infty}\max\{\nu_1(n),\ldots,\nu_m(n)\}=0$ by
Theorem \ref{Main Thm} (ii). In fact, note that if $\mathbf{A}_n$
and $K_{n;\nu_1,\ldots,\nu_m}$ have the same partition, $F_1\equiv
0$ and $F_2$ is a Bernoulli distribution with mean $p$, then
$\mathbf{A}_n$ is the adjacent matrix of
$G_{n;\nu_1\ldots\nu_m}(p)$. Thus, by Theorem \ref{Main Thm} (ii),
a.e. random graph $G_{n;\nu_1\ldots\nu_m}(p)$ enjoys the following
equation:
$$  \En(G_{n;\nu_1\ldots\nu_m}(p))
  =n^{3/2}\left(\frac{8}{3\pi}\sqrt{p(1-p)}+o(1)\right).
$$

For $m$-partite random graphs $G_{n;\nu_1\ldots\nu_m}(p)$ such that
$$\lim_{n\rightarrow\infty}\max\{\nu_1(n),\ldots,\nu_m(n)\}>0
\mbox{ and there exist $\nu_i$ and $\nu_j$ such that }
\lim_{n\rightarrow\infty}\frac{\nu_i(n)}{\nu_j(n)}<1,
$$
we can establish lower and upper bounds for its energy. For the
purpose, we first introduce the following an auxiliary assertion due
to \cite{ky}.
\begin{lem}\label{ky}
Let $\mathbf{X},\mathbf{Y},\mathbf{Z}$ be square matrices of order
$n$ such that $\mathbf{X}+\mathbf{Y}=\mathbf{Z}$, then
$$\sum_{i=1}^ns_i(\mathbf{X})+\sum_{i=1}^ns_i(\mathbf{Y})\geq\sum_{i=1}^ns_i(\mathbf{Z})$$
where $s_i$ $(i=1,\cdots,n)$ is the singular values of a matrix.
\end{lem}

Similarly, suppose $\mathbf{A}_n$ and $\mathcal
{G}_{n;\nu_1\ldots\nu_m}(p)$ have the same partition
$V_1,\cdots,V_m$ $(|V_i|=\nu_in)$. Then, $\mathbf{A}_n$ is the
adjacent matrix of $G_{n;\nu_1\ldots\nu_m}(p)$ providing $F_1\equiv
0$ and $F_2$ is a Bernoulli distribution $B(p)$.  Without loss of
generality, we assume, for some $r\ge 1$, $|V_1|,\ldots,|V_r|$ are
of order $O(n)$ while $|V_{r+1}|,\cdots,|V_m|$ of order $o(n)$. Let
$\mathbf{X}'_n$ be a random symmetric matrix such that

$$\mathbf{X}'_n(ij)=\left\{\begin{array}{ll}\mathbf{A}_n(ij) &
\mbox{if  $i$ or $j\notin V_s(1\leq s\leq r)$,}\\
B(p) & \mbox{if $i,j\in V_s(1\leq s\leq r)$ and $i>j$,}\\
0 & \mbox{if $i,j\in V_s(1\leq s\leq r)$ and $i=j.$}\\
\end{array}\right.$$   From  Theorem \ref{Main Thm} (ii),
$\mathbf{X}'_n$ has the same LSD as $\mathbf{X}_n$ on condition that
$F_1\equiv 0$ and $F_2=B(p)$.  Set
\begin{equation}\label{dxa}
\mathbf{D}_n=\mathbf{X}'_n-\mathbf{A}_n=
\left(\begin{array}{lllll} \mathbf{K}_1 & & & &\\
 & \mathbf{K}_2& & & \\
 & &\ddots& &\\
 & & & \mathbf{K}_r&\\
 & & & & \mathbf{O}
\end{array}\right)_{n\times n}\end{equation}
Let $\mathbf{M}$ be a matrix. We use $\En(\mathbf{M})$ to denote the
sum of singular values of $\mathbf{M}$. Evidently, if $\mathbf{M}$
is the adjacent matrix of a simple graph $G$ then
$\En(G)=\En(\mathbf{M}).$ One can readily see that a.e. matrix
$\mathbf{K}_i~(i=1,\ldots,r)$ enjoys the following:
$$\En(\mathbf{K}_i)=\left(\frac{8}{3\pi}\sqrt{p(1-p)}+o(1)\right)(\nu_in)^{3/2},$$
 and then a.e. matrix
$\mathbf{D}_n$ satisfies the following:
$$\En(\mathbf{D}_n)=\left(\frac{8}{3\pi}\sqrt{p(1-p)}+o(1)\right)\left(\nu_1^{\frac{3}{2}}
+\cdots+\nu_r^{\frac{3}{2}}\right)n^{\frac{3}{2}}.$$ By (\ref{dxa}),
we have $\mathbf{A}_n+\mathbf{D}_n=\mathbf{X}'_n$ and
$\mathbf{X}'_n+(-\mathbf{D}_n)=\mathbf{A}_n$. Employing Lemma
\ref{ky}, we deduce
$$\En(\mathbf{X}'_n)-\En(\mathbf{D}_n)\leq \En(\mathbf{A}_n)\leq \En(\mathbf{X}'_n)+\En(\mathbf{D}_n).$$
Therefore, we establish the following result.
\begin{thm}\label{budengbu}
Let $G_{n;\nu_1\ldots\nu_m}(p)$ be a random graph of
$\mathcal{G}_{n;\nu_1\ldots\nu_m}(p)$. Then a.e. random graph
$G_{n;\nu_1\ldots\nu_m}(p)$ satisfies the following inequality
$$\left(1-\sum_{i=1}^r\nu_i^{\frac{3}{2}}\right)n^{3/2}
\le
\En(G_{n;\nu_1\ldots\nu_m}(p))\left(\frac{8}{3\pi}\sqrt{p(1-p)}+o(1)\right)^{-1}
\le\left(1+\sum_{i=1}^r\nu_i^{\frac{3}{2}}\right)n^{3/2}
.$$
\end{thm}

\noindent\textbf{Remark}. Since $\nu_1,\ldots,\nu_r$ are positive
real numbers with $\sum_{i=1}^r\nu_i\le 1$, we have
$\sum_{i=1}^r\nu_i(1-\nu_i^{1/2})>0$. Therefore,
$\sum_{i=1}^r\nu_i>\sum_{i=1}^r\nu_i^{3/2},$ and thus
$1>\sum_{i=1}^r\nu_i^{3/2}$. Hence, we can deduce, by the theorem
above, that a.e. random graph $G_{n;\nu_1\ldots\nu_m}(p)$ enjoys the
following $$\En(G_{n;\nu_1\ldots\nu_m}(p))=O(n^{3/2}).$$

\end{document}